\RequirePackage{ifpdf}
\ifpdf 
\documentclass[pdftex]{sigma}
\else
\documentclass{sigma}
\fi

\usepackage[all]{xy}

  \def\<{{\langle}}
  \def\>{{\rangle}}

  \def\eps{\varepsilon}

  \def\note#1{{}}

  \def\note#1{}
  
  \def\N{\mathbb{N}}

  \def\cD{{\mathcal D}}

  \def\rhom#1#2#3{{{\rm Hom}\sb{#1}(#2,#3)}}

  \def\C{\mathbb{C}}
   \def\R{\mathbb{R}}
  \def\Z{\mathbb{Z}}

  \def\beq{\begin{equation}}
  \def\eeq{\end{equation}}

  \def\im{{\rm Im}}

  \def\ot{{\otimes}}

  \def\coker{\mathrm{coker}}

    \def\oan#1{\Omega^{#1}(A)}
    \def\oa{\Omega (A)}

     \def\1{\mathbf{1}}

\def\e2{E_q(2)}
\def\ome#1{\Omega^#1(\e2)}
\def\omc#1{\Omega^#1(\C_q)}
\def\integ#1#2{\mathfrak{I}^{#1}(#2)}

\def\k{\Bbbk}

  \newcounter{zlist}

  \newcounter{blist}

  \newcounter{rlist}

   \newcounter{alist}

\def\stac#1{\raise-.2cm\hbox{$\stackrel{\displaystyle\otimes}{\scriptscriptstyle{#1}}$}}

\def\cten#1{\raise-.2cm\hbox{$\stackrel{\displaystyle\widehat{\otimes}}
{\scriptscriptstyle{#1}}$}}

  \def\Label#1{\label{#1}\ifmmode\llap{[#1] }\else
  \marginpar{\smash{\hbox{\tiny [#1]}}}\fi}
  \def\Label{\label}

\begin{document}

\renewcommand{\thefootnote}{$\star$}

\renewcommand{\PaperNumber}{040}

\FirstPageHeading

\ShortArticleName{Integral calculus on $E_q(2)$}

\ArticleName{Integral calculus on $\boldsymbol{E_q(2)}$\footnote{This paper is a
contribution to the Special Issue ``Noncommutative Spaces and Fields''. The
full collection is available at
\href{http://www.emis.de/journals/SIGMA/noncommutative.html}{http://www.emis.de/journals/SIGMA/noncommutative.html}}}

\Author{Tomasz BRZEZI\'NSKI}

\AuthorNameForHeading{T.\ Brzezi\'nski}

\Address{Department of Mathematics, Swansea University,
  Singleton Park, Swansea SA2 8PP, UK}
\Email{\href{mailto:T.Brzezinski@swansea.ac.uk}{T.Brzezinski@swansea.ac.uk}}
\URLaddress{\url{http://www-maths.swan.ac.uk/staff/tb/}}

\ArticleDates{Received March 11, 2010;  Published online May 13, 2010}

\Abstract{The complexes of integral forms on the quantum Euclidean group $E_q(2)$ and the quantum plane are def\/ined and their isomorphisms with the corresponding de Rham complexes are established.}


\Classification{58B32; 16T05; 81R60; 81R50}

\section{Introduction}
In \cite[Chapter~4]{Man:gau} the Berezin integral on a supermanifold is explained in terms of the complex of {\em integral forms}. The boundary operator on this complex is derived from a~unique {\em right connection}, a notion originating from the analysis of  supersymmetric $\cD$-modules \cite{Pen:mod}. A right connection is a~version of a~covariant derivative, whose covariance properties are dif\/ferent from (one might say: dual to) that of the usual connection. The existence of a right connection can be proven by a direct construction in terms of  local coordinates.

The notions of a right connection and a complex of integral forms were extended to dif\/ferential operators over  commutative algebras in \cite{Vin:coc} and, further,
over a commutative algebra in the braided category of Yetter--Drinfeld modules in \cite{Ver:dif}.

Studies of right connections  for non-commutative algebras were initiated  in \cite{Brz:con}. Right connections  for a non-commutative algebra are def\/ined as maps between modules of (right) linear homomorphisms from (homogeneous parts of) a dif\/ferential graded algebra to a module over the zero degree part of this algebra, and hence termed  {\em hom-connections} (as opposed to usual connections which operate between the tensor product of a module and a dif\/ferential graded algebra). The local arguments which lead to the existence of right connections on supermanifolds are not available in the realm of noncommutative geometry. In their stead an ef\/fective method of constructing dif\/ferential graded algebras with a unique hom-connection was presented in~\cite{BrzElK:int}. This method is based on the use of twisted multi-derivations.

The aim of this note is to calculate explicitly two examples of complexes of integral forms based on twisted multi-derivations, and to show that these complexes are isomorphic to the corresponding de Rham complexes, which implies (a version of) the Poincar\'e duality. The examples in question are the quantum Euclidean group $\e2$ obtained by the contraction  of the quantum group $SU_q(2)$ in \cite{Wor:con}, and its homogeneous space -- the quantum plane.

The paper is organised as follows. In Section~\ref{sec.fund} we brief\/ly recall from \cite{Brz:con} and \cite{BrzElK:int} def\/initions of hom-connections and integral forms and the method of constructing hom-connections in dif\/ferential calculi determined by $q$-skew derivations. In Section~\ref{sec.e2} we recall from \cite{Wor:con} the def\/inition of the quantum group $\e2$, describe a three-dimensional dif\/ferential calculus on it and then explicitly construct a f\/lat hom-connection and the complex of integral forms on $\e2$. The explicit isomorphism of the complexes of dif\/ferential and integral forms is then presented. In Section~\ref{sec.cq} we derive a hom-connection and the complex of integral forms on the quantum plane understood as the (quantum) homogeneous space of $\e2$. The paper is concluded with comments and outlook.

We hope that the material contained in Section~\ref{sec.fund} of this note will provide the reader
with a~concise introduction to non-commutative integral forms.
We also hope that explicitly calculated examples (in Sections~\ref{sec.e2} and~\ref{sec.cq}),
which append examples presented in \cite{BrzElK:int},
will indicate the methods and techniques involved in construction of integral forms on quantum spaces.

\section{Noncommutative integral forms}\label{sec.fund}

An algebra means an associative algebra with identity over a f\/ield $\k$.
By a {\em differential graded algebra over an algebra $A$} we mean a non-negatively graded algebra
\begin{displaymath}
\oa := \bigoplus_{k=0}^\infty \oan k,
\end{displaymath}
such that $\oan 0 =A$, together with a degree one operation $d: \oan k\to \oan {k+1}$ that is nilpotent, i.e.\ $d\circ d = 0$, and satisf\/ies the graded Leibniz rule, i.e., for all $\omega \in \oan k$ and $\omega'\in \oa$,
\begin{displaymath}
d(\omega \omega ') = d\omega \omega' + (-1)^k \omega d\omega'.
\end{displaymath}
Each of the $\oan k$ is an $A$-bimodule, and we denote the right dual of $\oan k$ (i.e.\ the collection of all right $A$-linear maps $\oan k \to A$) by $\integ k A$, that is
\begin{displaymath}
\integ kA := \rhom A {\oan k} A.
\end{displaymath}
Each of the $\integ kA$ is an $A$-bimodule with the actions def\/ined by
\begin{equation}\label{def.act.integ}
a\cdot \phi \cdot b (\omega) := a\phi(b\omega), \qquad \mbox{for all} \  \ a,b\in A, \  \phi\in \integ kA, \  \omega \in \oan k.
\end{equation}
To ease the notation we write $a \phi b$ for $a\cdot \phi \cdot b$. The right action \eqref{def.act.integ} is a special case of more general operation that, for any $\omega\in \oan l$, sends $\phi\in \integ{k+l} A$ to $\phi\omega \in \integ{k} A$. The latter is def\/ined by
\begin{displaymath}
\phi\omega (\omega') := \phi (\omega\omega'), \qquad \mbox{for all} \  \ \omega' \in \oan k.
\end{displaymath}

A {\em hom-connection} on a $\k$-algebra $A$ (with respect to the dif\/ferential graded algebra~$\oa$ over~$A$) is a $\k$-linear map $\nabla: \integ 1 A \to A$ such that,
\begin{displaymath}
\nabla (\phi a) = \nabla(\phi) a + \phi(da), \qquad \mbox{for all} \  \ a\in A, \  \phi\in \integ 1A.
\end{displaymath}
Any hom-connection can be extended to a family of maps $\nabla_k : \integ {k+1} A \to \integ k A$ by the formula
\begin{displaymath}
\nabla_k(\phi)(\omega) :=  \nabla  (\phi\omega) + (-1)^{k+1} \phi(d\omega), \qquad \mbox{for all} \ \ \phi\in \integ {k+1} A, \ \omega \in \oan{k}.
\end{displaymath}
 The $\nabla_k$ satisfy the following graded Leibniz rule
 \begin{equation}\label{hom.Leibniz}
 \nabla_l (\phi \omega) = \nabla_{k+l}(\phi)\omega + (-1)^{k+l}\phi d\omega, \qquad \mbox{for all} \ \ \phi\in \integ {k+l+1} A, \ \omega \in \oan{k}.
 \end{equation}
 The map $F := \nabla \circ \nabla_1$ is a right $A$-module homomorphism which is called the {\em curvature} of~$\nabla$. A~hom-connection~$\nabla $ is said to be {\em flat} provided $F=0$. To a f\/lat hom-connection $\nabla $ one associates a chain complex $\left(\bigoplus_{k=0} \integ k A, \nabla\right)$. This complex is termed a {\em complex of integral forms} on~$A$, and the canonical map
 \begin{displaymath}
 \Lambda: \ A \longrightarrow  \coker (\nabla)= A/\im\nabla
 \end{displaymath}
 is called a {\em $\nabla$-integral on $A$} (or simply an integral on $A$).

A  general construction of hom-connections based on twisted multi-derivations was presented in \cite{BrzElK:int}. This construction is applicable to all left-covariant dif\/ferential calculi on quantum groups. Presently we outline this construction in the special case of {\em $q$-skew derivations} and the reader is referred to \cite[Section~3]{BrzElK:int} for more details and for the general case.

Let $A$ be a $\k$-algebra. Following  \cite{GooLet:pri}, a linear map $\partial : A\to A$ is called a {\em $q$-skew derivation} provided there exists an algebra automorphism $\sigma: A\to A$ and a scalar $q$ such that $\sigma^{-1} \circ \partial\circ \sigma = q\partial$ and, for all $a,b\in A$,
\begin{equation}\label{skew}
\partial (ab) = \partial(a) \sigma(b) + a\partial(b),
\end{equation}
i.e.\ $\partial$ satisf\/ies the {\em $\sigma$-twisted Leibniz rule}. Starting with $q_i$-skew derivations $\partial_i$, $i=1,2,\ldots, n$, (with corresponding automorphisms $\sigma_i$) one constructs a f\/irst order dif\/ferential calculus on $A$ as follows. $\oan 1$ is a free left $A$-module $\bigoplus_{i=1}^n A \omega_i$  with basis $\omega_1,\ldots, \omega_n$ and right $A$-action given~by
\begin{equation}\label{eq.dif.rel}
 \omega_i a =  \sigma_{i}(a) \omega_i ,  \qquad i=1,2,\ldots, n.
\end{equation}
The exterior dif\/ferential $d: A\to \oan 1$ is def\/ined by the formula
\begin{equation}\label{eq.dif.d}
da = \sum_i \partial_i(a)\omega_i = \sum_{i,j}\omega_i \sigma_{i}^{-1}\left(\partial_j\left(a\right)\right).
\end{equation}
Often one reserves the term {\em differential calculus} to the pair $(\oan 1, d)$ that is {\em dense} in the sense that every element of $\oan 1$ is of the form $\sum_i a_idb_i$, for some $a_i, b_i\in A$. The calculus described above is dense if and only if there exist two f\/inite subsets $\{a_{i\,t}\}$, $\{b_{i\,t}\}$ of elements of $A$  such that,
\begin{displaymath}
\sum_t a_{i\,t}\partial_k(b_{i\,t}) = \delta_{ik}, \qquad \text{for all}  \ \ i,k =1, \ldots, n.
\end{displaymath}
The examples of calculi discussed in Sections~\ref{sec.e2} and \ref{sec.cq} are dense in this sense.

$(\oan 1, d)$  can be extended to a full dif\/ferential graded algebra $(\oa ,d)$ in the standard way, that is by using the graded Leibniz rule and the $A$-bimodule structure of $\oan 1$.

The calculus $(\oan 1, d)$ determined by the formulae \eqref{eq.dif.rel}, \eqref{eq.dif.d} admits a hom-connection $\nabla : \integ 1 A \to A$, given by, for all $f\in \integ 1 A $,
\begin{equation}\label{nabla}
\nabla (f) = \sum_i q_i \partial_i (f(\omega_i)).
\end{equation}
This is a unique hom-connection on $\oa$ with the property that $\nabla(\xi_i) =0$, where
$\xi_i\in \integ 1 A$ are  such that $\xi_i(\omega_j) = \delta_{ij}$, $i,j=1,2,\ldots , n$.

\section[Integral geometry of $E_q(2)$]{Integral geometry of $\boldsymbol{E_q(2)}$}\label{sec.e2}

In this and the following section the algebras are over the f\/ield of complex numbers.
The (polynomial part of the) quantum Euclidean group $\e2$, obtained by the contraction of $SU_q(2)$ in \cite{Wor:con}, is a $*$-Hopf algebra generated by $v$, $n$ subject to relations
\begin{displaymath}
nn^* = n^*n, \qquad vv^* = v^*v=1, \qquad vn=qnv, \qquad vn^* = qn^*v,
\end{displaymath}
where $q\in \R$. The elements $v$, $v^*$ are grouplike, while
\begin{displaymath}
\Delta(n) = v\ot n + n\ot v^*, \qquad \eps(n) =0, \qquad S(n)= -q^{-1}n, \qquad S(n^*)= -qn^*.
\end{displaymath}
By setting $z = vn$, the algebra $\e2$ can be equivalently def\/ined as generated by a unitary element $v$ and elements $z$, $z^*$ subject to relations
\begin{equation}\label{eq2.cone}
vz = qzv, \qquad vz^* = qz^*v, \qquad  zz^* = q^2z^*z.
\end{equation}
The algebra $\e2$ has a $\Z$-grading, def\/ined by
\begin{displaymath}
|v| = |n^*|=1, \qquad |v^*| = |n|=-1  .
\end{displaymath}
In particular, the elements $z$, $z^*$ generate the zero-degree subalgebra of $\e2$ which we denote by $\C_{q}$.

A three-dimensional left-covariant calculus $\ome 1$ on $\e2$ can be obtained by contraction of the 3D calculus on $SU_q(2)$ introduced in~\cite{Wor:twi}. $\ome 1$ is generated by left-invariant forms~$\omega_0$,~$\omega_\pm$ subject to relations
\begin{alignat}{5}
& \omega_0  v =q^{-2}v   \omega_0 , \qquad&&
\omega_0  n= q^2n   \omega_0, \qquad&&
\omega_0  n^* =q^{-2}n^*  \omega_0 , \qquad&&
\omega_0  v^*= q^2v^*  \omega_0 , &  \nonumber\\
& \omega_\pm   v =q^{-1}v  \omega_\pm , \qquad&&
\omega_\pm   n= qn  \omega_\pm, \qquad&&
\omega_\pm   n^* =q^{-1}n^* \omega_\pm , \qquad&&
\omega_\pm   v^*= qv^*   \omega_\pm .&\label{eq.3Drel}
\end{alignat}
This is a $*$-calculus with $\omega_0^* = -\omega_0$, $\omega_\pm^* = q^{\mp 1}\omega_\mp$. The action of the exterior dif\/ferential $d$ on the generators
is
\begin{alignat}{3}
& d v = v \omega_0,\qquad&&
d n = -q^2 n \omega_0 +v \omega_- ,& \nonumber\\
& d n^* = n^*\omega_0 + q^2 v^*\omega_+ ,\qquad&&
d v^* = -q^2 v^*\omega_0 .& \label{eq.3Drel.d}
\end{alignat}
Consequently,
\begin{equation}\label{eq.omega}
\omega_0 = v^* dv, \qquad \omega_- = v^* dn - q^{-1}n dv^*, \qquad \omega_+ = q^{-2} v dn^* - q^{-1} n^* dv.
\end{equation}

In view of commutation rules \eqref{eq.3Drel}, the module structure of $\ome 1$ is of the type described by \eqref{eq.dif.rel} with automorphisms $\sigma_0$ and $\sigma_+=\sigma_-$ given by, for all homogeneous $a\in A$ (with the $\Z$-degree $|a|$),
\begin{equation}\label{sigma.3d}
\sigma_0(a) = q^{-2|a|}a, \qquad \sigma_\pm(a) = q^{-|a|}a.
\end{equation}
The formulae \eqref{eq.3Drel.d} or \eqref{eq.omega} indicate that the $\Z$-grading of $\e2$ can be extended to a $\Z$-grading of $\ome 1$ such that the dif\/ferential $d$ is the degree preserving map, by setting
\begin{displaymath}
|\omega_0| =0, \qquad |\omega_\pm | = \pm 2.
\end{displaymath}
Equations \eqref{eq.3Drel.d} determine (or can be understood as determined by as in \eqref{eq.dif.d})  maps $\partial_i: \e2\to\e2$ that satisfy the $\sigma_i$-twisted derivation properties \eqref{skew}. Explicitly, in terms of actions on generators of $\e2$ the maps $\partial_i$ are
\begin{alignat*}{5}
& \partial_0(v) = v,  \qquad&&   \partial_0(n) = -q^2v,  \qquad&&  \partial_0(n^*) = n^*,   \qquad && \partial_0(v^*) = -q^2v^*  ,& \\
& \partial_+(v) = 0, \qquad&&  \partial_+(n) = 0,  \qquad&&  \partial_+(n^*) = q^2v^*,  \qquad && \partial_+(v^*) = 0  ,& \\
& \partial_-(v) = 0, \qquad&&  \partial_-(n) = v, \qquad && \partial_-(n^*) = 0, \qquad  && \partial_-(v^*) = 0  .&
\end{alignat*}
Since $da = \partial_-(a)\omega_- + \partial_0(a)\omega_0 + \partial_+(a)\omega_+$, the derivations have $\Z$-degrees
\begin{equation}\label{partial.deg}
 |\partial_0| =0, \qquad |\partial_\pm | = \mp 2,
\end{equation}
i.e., for all homogeneous $a\in \e2$, $|\partial_0(a)| = |a|$, $|\partial_\pm(a)| = |a | \mp 2$. Combining equations \eqref{sigma.3d} with \eqref{partial.deg} one easily f\/inds that the maps  $\partial_0$, $\partial_+$, $\partial_-$ are $q_i$-skew derivations with constants~1, $q^{-2}$ and $q^2$, respectively. Therefore, there is a hom-connection on~$\e2$ (with respect to~$\ome 1$) given by the formula~\eqref{nabla}, that is, for all $f\in \integ 1 {\e2}$,
\begin{equation}\label{eq.hom.3d}
\nabla  (f) = \partial_0\left( f\left(\omega_0\right)\right) + q^{-2}\partial_+\left( f\left(\omega_+\right)\right) + q^{2}\partial_-\left( f\left(\omega_-\right)\right).
\end{equation}

The calculus $\ome 1$ can be extended to the full dif\/ferential graded algebra. The relations in the bimodule $\ome 2$ (deformed exterior product) are
\begin{equation}\label{eq.3Drel.higher}
\omega_i^2  =  0, \qquad
 \omega_+\omega_-=-q^2\omega_-\omega_+ ,\qquad
\omega_0\omega_-=-q^4\omega_-\omega_0 ,\qquad
\omega_+\omega_0=-q^4\omega_0\omega_+ ,
\end{equation}
and the exterior derivative is
\begin{equation}\label{eq.3Drel.d.higher}
d \omega_0 =0 ,\qquad
d \omega_+ =q^2(q^2+1)\,\omega_0\omega_+ , \qquad
d \omega_- =q^2(q^2+1)\omega_-\omega_0 .
\end{equation}
Thus $\ome 2$ is a free module generated by three closed forms $\omega_-\omega_0$, $\omega_0\omega_+$ and $\omega_-\omega_+$. In degree $3$, $\ome 3$ is generated by the  form $\omega_-\omega_0\omega_+$.

It might be worth noticing at this point that the dif\/ferential graded algebra $\ome\bullet$ can be equivalently (and conveniently) described in terms of generators $v$, $z$, $z^*$. An easy calculation reveals that
\begin{equation}\label{omega.z}
\omega_+=q^{-3} v^2dz^*, \qquad \omega_- = {v^*}^2dv,
\end{equation}
so that
\begin{equation}\label{omega.dz}
dz^* = q^3 {v^*}^2 \omega_+, \qquad dz = v^2\omega_-.
\end{equation}
Therefore, the $*$-calculus $\ome 1$ is freely generated by $dv$, $dz$ and $dz^*$ subject to relations
\begin{alignat}{5}
& vdv = q^2dv v, \qquad && v^*dv = q^{-2}dvv^*, \qquad && zdv = q^{-1}dv z, \qquad && z^*dv = q^{-1}dv z^*, & \nonumber\\
& vdz = qdz v, \qquad && vdz^* = qdz^*v, \qquad && zdz = q^{-2}dz z, \qquad && z^*dz = q^{-2} dzz^*. & \label{rel.z}
\end{alignat}
The wedge product calculated from relations \eqref{eq.3Drel.higher} comes out as
\begin{alignat}{3}
& (dz)^2   =  (dv)^2   =  dvdv^*    =   0, \qquad  && dvdz = -qdzdv,&   \nonumber\\
 & dvdz^*=-q dz^*dv,  \qquad  && dzdz^* = - q^2 dz^*dz.& \label{wedge.z}
\end{alignat}
Note that relations \eqref{wedge.z} are simply obtained by dif\/ferentiating rules \eqref{rel.z} (and observing that $v^*$ is the inverse of $v$).

Every element of $\ome 3$ can be written as a linear combination of
\begin{displaymath}
v^k{z^*}^{l-1}z^{m-1} dzdvdz^*, \qquad k\in \Z, \  l,m\in \N .
\end{displaymath}
On the other hand, relations  \eqref{rel.z} and the Leibniz rule imply that
\begin{displaymath}
d\big(v^k{z^*}^{l-1}z^{m}\big) = d\big(v^k{z^*}^{l-1}\big)z^{m} + [m]_{q^2}v^k{z^*}^{l-1}z^{m-1}dz,
\end{displaymath}
where
\begin{displaymath}
[m]_x = 1+ x +\dots + x^{m-1}
\end{displaymath}
is a notation for the $x$-integer. Therefore,{\samepage
\begin{equation}\label{el.hom}
v^k{z^*}^{l-1}z^{m-1} dzdvdz^* = d\left(\frac{1}{[m]_{q^2}} v^k{z^*}^{l-1}z^{m} dvdz^*\right),
\end{equation}
so the third de Rham cohomology group of $\e2$ is trivial, $H^3(\e2)=0$.}

Our next aim is to show that $\nabla$ def\/ined by equation \eqref{eq.hom.3d} is a f\/lat hom-connection and thus def\/ines a complex of integral forms on $\e2$. To this end we need f\/irst to extend $\nabla$ to $\nabla_1: \integ 2 {\e2}  \to \integ 1 {\e2} $. Def\/ine  $\phi_0, \phi_\pm \in \integ 2 {\e2}$ by
\begin{displaymath}
\phi_0(\omega_-\omega_+) =1  ,  \qquad \phi_+(\omega_-\omega_0) =1  , \qquad   \phi_-(\omega_0\omega_+) =1  ,
\end{displaymath}
and zero on other generators. By inspection of relations~\eqref{eq.3Drel.higher} one concludes that any $f\in \integ 2 {\e2}$ can be written as
$
f = \phi_{0}a_0 + \phi_+a_+ +\phi_- a_-,$ for suitably def\/ined $a_i\in \e2$. Since the curvature of a hom-connection is a right $\e2$-linear map, it suf\/f\/ices to compute it on the~$\phi_{i}$. Let $\xi_0, \xi_\pm \in \integ 1 {\e2}$ be the right dual basis to $\omega_0$, $\omega_\pm$, i.e.\ such that $\xi_0(\omega_0) = \xi_\pm (\omega_\pm) = 1$ and zero on other generators of $\ome 1$. Using commutation rules \eqref{eq.3Drel.higher} one easily computes that
\begin{alignat}{4}
&\phi_0\omega_0 =0, \qquad&& \phi_0\omega_+ = - q^2\xi_-, \qquad &&\phi_0\omega_- = \xi_+, &\nonumber\\
&\phi_+\omega_0 =-q^4\xi_-, \qquad&& \phi_+\omega_+ = 0, \qquad &&\phi_+\omega_- = \xi_0,& \nonumber\\
&\phi_-\omega_0 =\xi_+, \qquad &&\phi_-\omega_+ = - q^4\xi_0, \qquad &&\phi_-\omega_- = 0.&\label{phi.omega}
\end{alignat}
By def\/inition, for all one-forms $\omega$ and all $\phi\in \integ 2 {\e2}$,
\begin{equation}\label{nab1.0}
\nabla_1(\phi)(\omega) = \nabla(\phi\omega) +\phi(d\omega).
\end{equation}
Combining \eqref{phi.omega} with \eqref{nab1.0} and \eqref{eq.3Drel.d.higher}, and remembering that $\nabla $ given by \eqref{eq.hom.3d} is the unique hom-connection such that $\nabla (\xi_i) =0$, $i=0,+,-$, we obtain
\begin{displaymath}
\nabla_1(\phi_0) =0, \qquad \nabla_1(\phi_+) = q^2(q^2+1)\xi_-, \qquad \nabla_1(\phi_-) = q^2(q^2+1)\xi_+.
\end{displaymath}
Since $\nabla (\xi_i) =0$,  $\nabla\circ\nabla_1(\phi_i) =0$, for all $i=0,+,-$, and the hom-connection \eqref{eq.hom.3d} is f\/lat.

The module $\integ 3 {\e2}$ is generated by $\phi\in \integ 3 {\e2}$ def\/ined by $\phi(\omega_-\omega_0\omega_+)=1$. Using commutation rules \eqref{eq.3Drel.higher}, one f\/inds
\begin{equation}\label{eq.psi}
\phi \omega_-\omega_0 = \xi_+, \qquad \phi \omega_0\omega_+ = q^6\xi_-, \qquad \phi \omega_-\omega_+ = -q^4\xi_0.
\end{equation}
Since, for any 2-form $\omega$,  $\nabla_2 (\phi)(\omega) =\nabla(\phi \omega) -\phi (d\omega)$, the forms generating $\ome 2$ are closed and $\nabla (\xi_i) =0$, $i=0,+,-$, we conclude that
\begin{displaymath}
\nabla_2(\phi) =0.
\end{displaymath}
This completes the description of the complex
\begin{equation}\label{comp.int}
\xymatrix{ 0 \ar[r] & \integ 3{\e2} \ar[r]^-{\nabla_2} &  \integ 2{\e2} \ar[r]^-{\nabla_1} & \integ 1{\e2} \ar[r]^-{\nabla} & {\e2} \ar[r]^-{\Lambda} & \coker{\nabla} \ar[r] & 0}
\end{equation}
of integral forms on $\e2$. Here $\Lambda$ is the associated integral on $\e2$. We will presently show that $\e2$ enjoys the strong Poincar\'e duality in the sense that  the complex~\eqref{comp.int} is isomorphic to the de Rham complex. In view of the triviality of $H^2(\e2)$ this will allow us to deduce that $\coker {\nabla} =0$, hence the integral $\Lambda$ is zero. The isomorphism of the de Rham and integral complexes means the following commutative diagram
\begin{equation}\label{diag.Poin.ex}
\xy<0mm,8mm> \xymatrix{ \e2 \ar[rr]^d \ar[d]_{\Theta^*} && \ome 1\ar[rr]^d \ar[d]_\Phi && \ome 2\ar[rr]^d \ar[d]_\Psi &&\ome 3\\
\integ 3{\e2} \ar[rr]^{\nabla_2} &&\integ 2{\e2}\ar[rr]^{\nabla_1} && \integ 1{\e2}\ar[rr]^\nabla && \e2 \ar[u]_\Theta \, ,} \endxy
\end{equation}
in which all columns are (right $A$-module) isomorphisms. These are def\/ined as follows. $\Theta^*(a) = \phi a$, $\Theta (a) = \omega_-\omega_0\omega_+a$, and
\begin{gather*}
\Phi(\omega_-a + \omega_0 b +\omega_+ c) = \phi_- a -q^4\phi_0b + q^6\phi_+c  ,
\\
\Psi(\omega_-\omega_0a + \omega_-\omega_+ b +\omega_0\omega_+ c) = \xi_+ a -q^4\xi_0b + q^6\xi_-c ,
\end{gather*}
for all $a,b,c\in \e2$ (compare the form of $\Psi$ with relations \eqref{eq.psi}). The commutativity of the diagram \eqref{diag.Poin.ex} can be checked by a straightforward albeit lengthy calculation.  For example, the commutativity of the leftmost square can be verif\/ied as follows. Take a homogeneous  $a\in \e2$ with the $\Z$-degree $|a|$. In view of the form of commutation rules \eqref{eq.3Drel} or the automorphisms~$\sigma_\pm$,~$\sigma_0$ in \eqref{sigma.3d}, and the grading of derivations $\partial_\pm$, $\partial_0$, one easily computes
\begin{gather*}
\nabla_2\circ \Theta^* (a) (\omega_-\omega_0)  = \nabla_2(\phi)(\omega_-\omega_0)a + \phi da (\omega_-\omega_0) \\
\phantom{\nabla_2\circ \Theta^* (a) (\omega_-\omega_0) }{} = \phi(\partial_+(a)\omega_+\omega_-\omega_0) = q^{4|a| -8}\phi(\omega_+\omega_-\omega_0)\partial_+(a) = q^{4|a| -2}\partial_+(a),
\end{gather*}
where the f\/irst equality follows by the Leibniz rule for a hom-connection \eqref{hom.Leibniz}, and the last two equalities by the right linearity and the def\/inition of $\phi$ combined with the  commutation rules~\eqref{eq.3Drel.higher}. On the other hand
\begin{gather*}
\Phi(da) (\omega_-\omega_0)  = \Phi(\partial_-(a)\omega_- + \partial_0(a)\omega_0 +\partial_+(a)\omega_+)(\omega_-\omega_0) \\
\phantom{\Phi(da) (\omega_-\omega_0)}{} = q^{|a|+4}\phi_+(\partial_+(a) \omega_-\omega_0) =q^{4|a| -2}\partial_+(a),
\end{gather*}
where the equations \eqref{sigma.3d}, the def\/initions of $\Phi$ and of $\phi_\pm$, $\phi_0$ were used. In the same way one computes
\begin{displaymath}
\nabla_2\circ \Theta^* (a) (\omega_-\omega_+) = -q^{4|a| +4}\partial_0(a)= \Phi(da) (\omega_-\omega_+)
\end{displaymath}
and
\begin{displaymath}
\nabla_2\circ \Theta^* (a) (\omega_0\omega_+) = q^{4|a| +8}\partial_-(a)= \Phi(da) (\omega_0\omega_+) .
\end{displaymath}
Therefore,
\begin{displaymath}
\nabla_2\circ \Theta^* =  \Phi\circ d,
\end{displaymath}
as required. The commutativity of the other squares in  diagram \eqref{diag.Poin.ex} is proven similarly.

The commutativity of \eqref{diag.Poin.ex} implies that $\coker \nabla$ is isomorphic to the third de Rham cohomo\-lo\-gy group $H^3(\e2)$, and thus allows one to compute the form of the former (and hence  the integral on $\e2$). Recall that $H^3(\e2)=0$, so also $\coker \nabla =0$, i.e.\ $\nabla$ is surjective, and the corresponding integral $\Lambda$ is zero. More specif\/ically, take $a\in \e2$. If $\omega\in \ome 2$ is such that $d\omega = \Theta (a)$ (and $\omega$ exists by the triviality of $H^3(\e2)$), then
\begin{displaymath}
\nabla(\Psi(\omega)) =a.
\end{displaymath}
For an element $a=v^k{z^*}^{l}z^{m}$ of a basis for $\e2$,
\begin{displaymath}
\Theta \big(v^k{z^*}^{l}z^{m}\big) = q^{-4k-8+m+l} v^{k-1}{z^*}^{l}z^{m}dzdvdz^*,
\end{displaymath}
where the def\/inition of $\Theta$ and relations \eqref{omega.z}, \eqref{rel.z} and  \eqref{eq2.cone} were used. In view of \eqref{el.hom} followed by \eqref{rel.z} and   \eqref{omega.dz} the corresponding $\omega$ comes out as
\begin{displaymath}
\omega=\frac{q^{-4k-8+m+l}}{[m+1]_{q^2}} v^{k-1}{z^*}^{l}z^{m+1} dvdz^* = \frac{q^{-k-6+2m+2l}}{[m+1]_{q^2}}\omega_0\omega_+ v^{k-2}{z^*}^{l}z^{m+1} ,
\end{displaymath}
therefore,{\samepage
\begin{displaymath}
v^k{z^*}^{l}z^{m} = \nabla\left(\frac{q^{-k+2m+2l}}{[m+1]_{q^2}}\xi_-  v^{k-2}{z^*}^{l}z^{m+1} \right) ,
\end{displaymath}
which explicitly proves that $\nabla$ is onto.}

Note that although the isomorphisms $\Phi$, $\Psi$, $\Theta$ and $\Theta^*$ are right $\e2$-linear, they are not $\e2$-bimodule maps, when the $\integ k {\e2}$ are viewed as $\e2$-bimodules by \eqref{def.act.integ}. However, being right $\e2$-module isomorphisms, they can be forced to produce new left $\e2$-actions on the  $\integ k {\e2}$  and also to make the complex of integral forms a dif\/ferential graded algebra, isomorphic to $\Omega (\e2)$. The basic integral forms $\phi$, $\phi_\pm$, $\phi_0$, $\xi_\pm$, $\xi_0$ themselves form a seven-dimensional {\em skeletal algebra of integral forms} on $\e2$. The multiplication table can be easily worked out from the def\/initions of isomorphisms  $\Phi$, $\Psi$, $\Theta$ and $\Theta^*$, and from the multiplication rules in $\Omega (\e2)$:
\begin{center}
\begin{tabular}{c|c|c|c|c|c|c|c|}
&$\phi$ & $\phi_-$ &$\phi_0$ & $\phi_+$ & $\xi_-$ &$\xi_0$ & $\xi_+$\\ \hline
$\phi$ & $\phi$ & $\phi_-$ & $\phi_0$ & $\phi_+$ & $\xi_-$ &$\xi_0$ & $\xi_+$\\ \hline $\phi_-$ & $\phi_-$ & 0 & $-q^{-4} \xi_+$ & $-q^{-2} \xi_0$ & $q^{-6} \phi$ & 0 & 0 \\ \hline
$\phi_0$ & $\phi_0$ & $\xi_+$ & 0 & $-q^{-4} \xi_-$ & 0 &$-q^{-4} \phi$  & 0 \\ \hline
$\phi_+$ & $\phi_+$ & $\xi_0$ & $\xi_-$ & 0 & 0 & 0  &$ \phi$ \\ \hline
$\xi_-$ & $\xi_-$ &$\phi$ &  0 & 0  & 0 & 0 & 0 \\ \hline
$\xi_0$ & $\xi_0$ & 0 & $-q^{-4}\phi$  & 0  & 0 & 0 & 0 \\ \hline
$\xi_+$ & $\xi_+$ & 0 & 0 &  $q^{-6}\phi$   & 0 & 0 & 0 \\ \hline
\end{tabular}
\end{center}

\section[Integral geometry of $\C_q$]{Integral geometry of $\boldsymbol{\C_q}$}\label{sec.cq}

The subalgebra $\C_q$ of $\e2$ generated by $z$, $z^*$ is a quantum homogeneous space of $\e2$. It is also a base algebra for the quantum principal bundle with the quantum total space $\e2$ and the structure group $U(1)$ identif\/ied with the Hopf algebra of Laurent polynomials in variable $u$, $\C[u,u^{-1}]$ ($u$ is a grouplike element). In other words, $\e2$ is a principal comodule algebra. The coaction of $\C[u,u^{-1}]$ on $\e2$ is given by the $\Z$-grading, i.e.\ for any $a\in \e2$ of $\Z$-degree $|a|$,
\begin{equation}\label{coaction}
a\longmapsto a\ot u^{|a|}.
\end{equation}
$\C[u,u^{-1}]$ has a natural $*$-Hopf algebra structure given by $u^* = u^{-1}$. With respect to this, the coaction \eqref{coaction} is a $*$-algebra map.

The calculus  described by equations \eqref{rel.z} and \eqref{wedge.z} restricts to the calculus on $\C_q$. By \cite[Theorem~4.3]{BrzElK:int}, the hom-connection \eqref{eq.hom.3d} restricts  to the hom-connection on the quantum plane~$\C_q$,
\begin{equation}\label{hom.plane}
\nabla: \ \ \integ 1 {\C_q} \to \C_q, \qquad \nabla{f} = q^2 \partial_-(\hat{f}(\omega_-)) + q^{-1}\partial_+ (\hat{f}(\omega_+)),
\end{equation}
where
\begin{displaymath}
\hat{f}(\omega_-) = f\big(\omega_-v^2\big){v^*}^2, \qquad \hat{f}(\omega_+) = f\big(\omega_+ {v^*}^2\big)v^2.
\end{displaymath}

The skew derivations $\partial_+$ and $\partial_-$ are not def\/ined on $\C_q$ but on $\e2$, hence \eqref{hom.plane} is not the optimal description of the hom-connection $\nabla$ on $\C_q$ as it seems to depend on the embedding of $\C_q$ into $\e2$. The calculus $\omc 1$ is freely generated by the holomorphic form $dz$ and the antiholomorphic form $dz^*$. Using the right $\C_q$-linearity of $f \in \integ 1 {\C_q}$ and commutation rules~\eqref{rel.z} one easily f\/inds that
\begin{equation}\label{hatf.z}
\hat{f}(\omega_-) = q^{-2}f(dz){v^*}^2, \qquad \hat{f}(\omega_+) = q^{-1}f(dz^*)v^2.
\end{equation}
Next, introduce (twisted) derivatives $\partial$, $\bar\partial$ associated to $dz$, $dz^*$ by the formula
\begin{displaymath}
d(a) = \partial(a)dz + \bar\partial(a)dz^*, \qquad \mbox{for all} \ \  a\in \C_q.
\end{displaymath}
Computing $d(a{v^*}^2)$ and $d(av^2)$ with the help of equations \eqref{omega.dz} and commutation rules in $\ome 1$, one easily f\/inds that, for all $a\in \C_q$,
\begin{displaymath}
d(a{v^*}^2) = q^2\partial(a)\omega_- +\bar\partial(a) dz^*{v^*}^2, \qquad d(av^2) = \partial(a)dz v^2 + q\bar\partial(a)\omega_+.
\end{displaymath}
This implies that
\begin{equation}\label{part.part}
\partial_- (a{v^*}^2) = q^2\partial(a), \qquad \partial_+(av^2) = q\bar\partial(a) \omega_+, \qquad  \mbox{for all}  \ \ a\in \C_q.
\end{equation}
Putting equations \eqref{hatf.z}, \eqref{part.part} and \eqref{hom.plane} together we obtain the following formula for a hom-connection on $\C_q$ given in terms of (anti)holomorphic forms,
\begin{equation}\label{hom.z.z}
\nabla(f) = q^2 \partial\left( f\left(dz\right)\right) + q^{-2} \bar\partial\left( f\left(dz^*\right)\right), \qquad \mbox{for all} \ \  f\in \integ 1{\C_q}.
\end{equation}

Working entirely in terms of $dz$ and $dz^*$,  the module of integral forms $\integ 1 {\C_q}$ is generated by~$\xi$,~$\bar\xi$ def\/ined by
\begin{displaymath}
\xi(dz) = \bar\xi(dz^*) = 1, \qquad \xi(dz^*) = \bar\xi(dz) = 0.
\end{displaymath}
The formula \eqref{hom.z.z} immediately implies that $\nabla(\xi) = \nabla(\bar\xi) = 0$.

The degree-two integral forms $\integ 2 {\C_q}$ are generated by $\psi$ determined by $\psi(dzdz^*) =1$, as one easily f\/inds that $\Omega^2(\C_q)$ is generated by $dzdz^*$; see relations \eqref{wedge.z}. The latter imply that
\begin{displaymath}
\psi dz = \bar\xi, \qquad \psi dz^* = -q^{-2}\xi,
\end{displaymath}
and thus arguments analogous to those in the case of $\e2$ af\/f\/irm f\/latness of $\nabla$.

Similarly to the quantum Euclidean group, the quantum plane $\C_q$ enjoys the strong Poincar\'e duality, i.e.\ there is a commutative diagram
\begin{displaymath}\label{diag.Poin.cq}
\xy<0mm,8mm> \xymatrix{ \C_q \ar[rr]^d \ar[d]_{\Theta^*} && \Omega^ 1(\C_q)\ar[rr]^d \ar[d]_\Phi && \Omega^ 2(\C_q)\\
\integ 2{\C_q} \ar[rr]^{\nabla_1} &&\integ 1{\C_q}\ar[rr]^{\nabla}  && \C_q \ar[u]_\Theta \, ,} \endxy
\end{displaymath}
in which vertical maps are right $\C_q$-module isomorphisms. These are given by, for all $a,b\in \C_q$,
\begin{displaymath}
\Theta(a) = dzdz^* a, \qquad \Phi(dz a+ dz^* a) = -\bar\xi a + q^{-2}\xi b, \qquad \Theta^*(a) = \psi a.
\end{displaymath}
The second de Rham cohomology group is trivial, hence, consequently, the integral on $\C_q$ is zero.

\section{Comments}

In this note we presented two explicit examples of complexes of integral forms. In relation to such specif\/ic examples one can pose two more general questions. First, the quantum group~$\e2$ and its dif\/ferential calculus are obtained by contraction of $SU_q(2)$ and a calculus on it. It might be interesting to study in general the ef\/fect of contraction on integral geometry or to develop the contraction procedure for integral forms. Second, together with examples in~\cite{BrzElK:int}, the examples presented here provide an indication that  the (strong) Poincar\'e duality in noncommutative geometry could be understood as  existence of an isomorphism between dif\/ferential and integral complexes. One could even venture a suggestion  that a (compact without bounda\-ry) noncommutative dif\/ferentiable or smooth manifold should be understood as a dif\/ferential graded algeb\-ra~$\oa$ over $A$ which is isomorphic (as a complex) to the complex of integral forms $(\integ \bullet A,\nabla)$ (with respect to~$\oa$).  It would be rather interesting to investigate, what classes of dif\/ferential graded algebras can be characterised by this property, and how this viewpoint  on the Poincar\'e duality compares  with the  (algebraic, homological) noncommutative Poincar\'e duality of Van den Bergh~\cite{Van:rel}; see also~\cite{Kra:Poi}.

\pdfbookmark[1]{References}{ref}
\LastPageEnding

\end{document}